\documentclass[9pt]{article}

\usepackage{amssymb, amsmath, latexsym, amscd}

\def\sref#1{Section~\ref{#1}}

\def\rref#1{Remark~\ref{#1}}
\def\cref#1{Corollary~\ref{#1}}
\def\pref#1{Proposition~\ref{#1}}
\def\lref#1{Lemma~\ref{#1}}
\def\eref#1{Example~\ref{#1}}

\def\aref#1{Assumption~\ref{#1}}
\def\cdref#1{Condition~\ref{#1}}

\newtheorem{thm}{Theorem}[section] 
\newtheorem{dfn}[thm]{Definition} 
 
\newtheorem{rmk}[thm]{Remark}

\newtheorem{prop}[thm]{Proposition} 
\newtheorem{lem}[thm]{Lemma}
\newtheorem{exs}[thm]{Examples}
\newtheorem{ex}[thm]{Example} 
\newtheorem{subex}[thm]{Subexample}

\newtheorem{ass}[thm]{Assumption} 
\newtheorem{cond}[thm]{Condition} 
\newtheorem{ntt}[thm]{Notation}

\def\Hm#1{H^{#1}_{A_{\widetilde{\Omega}}} (f_{\Omega}^{-1}(-\infty, 
c + \delta))}   
\def\Hmpr#1{H^{#1}_{A_{\widetilde{\Omega}'}} (f_{\Omega'}^{-1}(-\infty, 
c' + \delta))}   
\def\Hp#1{H^{#1}_{A_{\widetilde{\Omega}}} (f_{\Omega}^{-1}(c - \delta , 
+\infty))}     
\def\Hppr#1{H^{#1}_{A_{\widetilde{\Omega}'}} (f_{\Omega'}^{-1}(c' - 
\delta , +\infty))}   
\def\Hpm#1{H^{#1}_{A_{\widetilde{\Omega}}} (f_{\Omega}^{-1}(c - \delta , 
c + \delta))}    
\def\Hpmpr#1{H^{#1}_{A_{\widetilde{\Omega}'}} (f_{\Omega'}^{-1}(c' - \delta , 
c' + \delta))}   
   
\def\Hpun#1{H^{#1}_{A_{\widetilde{\Omega  \cup \Omega'}}} (f_{\Omega  
\cup \Omega'}^{-1}( \tilde{c} - \delta , +\infty))}     
\def\Hpmun#1{H^{#1}_{A_{\widetilde{\Omega  \cup \Omega'}}} (f_{\Omega  
\cup \Omega'}^{-1}(\tilde{c} - \delta , \tilde{c} + \delta))}    
\def\Ho#1#2{H^{#1}_{A_{\widetilde{\Omega}}} (#2)}   
\def\Hopr#1#2{H^{#1}_{A_{\widetilde{\Omega}'}} (#2)}   
\def\Houn#1#2{H^{#1}_{A_{\widetilde{\Omega \cup \Omega'}}} (#2)}

\def\bi#1{b^{#1}_{A, \Omega} (c, \delta)}   
\def\bipr#1{b^{#1}_{A, \Omega'} (c', \delta)}   
\def\biun#1{b^{#1}_{A, \Omega \cup \Omega'} (\alpha c + (1 - \alpha) 
c', \delta)}   
\def\bil#1#2{b^{#1}_{A, \Omega_{#2}} (c, \delta)}   

\def\phi{\varphi_{A, \Omega, c, \delta}}

\def\ph#1{\varphi^{#1}_{A, \Omega, c, \delta}}   
\def\phpr#1{\varphi^{#1}_{A, \Omega', c', \delta}}   
\def\phun#1{\varphi^{#1}_{A, \Omega \cup \Omega', \alpha c + (1- \alpha) c', \delta}}   

\DeclareMathOperator{\rk}{rank}
\DeclareMathOperator{\Hom}{Hom}
\DeclareMathOperator{\supp}{supp}
\DeclareMathOperator{\Int}{int}
\DeclareMathOperator{\Ker}{Ker}
\DeclareMathOperator{\diam}{diam}
\DeclareMathOperator{\Crit}{Crit}
\DeclareMathOperator{\Dim}{dim}
\DeclareMathOperator{\card}{card}
\DeclareMathOperator{\cri}{cri}

\newcommand{\Pf}{\noindent{\bf Proof}. }

\newcommand{\cqfd}
{%
\mbox{}%
\nolinebreak%
\hfill%
\rule{2mm}{2mm}%
\medbreak%
\par%
}
\renewcommand{\a}{{\cal A}{}}

\newcommand{\C}{\mathbb C}

\newcommand{\R}{\mathbb R}

\newcommand{\T}{{\cal T}{}}

\newcommand{\Z}{\mathbb Z}
\newcommand{\eps}{\varepsilon }
\newcommand{\ol}{\overline}

\newcommand{\wrt}{with respect to }
\newcommand{\rp}{respectively }

\begin{document} 

\title{Dynamical Morse entropy}

\author{M\'elanie Bertelson$^\ast$\footnote{Research supported by a FNRS 
Charg\'ee de Recherches contract } \& Misha Gromov$^\dagger$} 

\date{April 30, 2004}

\maketitle

\begin{center}
$^\ast$ D\'epartement de Math\'ematiques \\
UCL\\
2, Chemin du Cyclotron \\
1348 Louvain-la-Neuve \\
Belgique\\
{\tt bertelson@math.ucl.ac.be}
\end{center}

\begin{center}
$^\dagger$ IH\'ES\\
35, route de Chartres \\
91440 Bures-sur-Yvette\\
France\\
{\tt gromov@ihes.fr}
\end{center}

\section{Introduction}

Consider a "crystal", that is, the standard lattice $\Gamma = \Z^3$ in $\R^3$ 
with identical "molecules" positioned at all sites (points) $\gamma$ in $\Gamma$. 
Denote by $M$ the configuration space of such a molecule which is assumed to be a 
smooth finite dimensional manifold and let $X = M^\Gamma$ be the configuration 
space of the crystal, that is, the infinite product of $\Gamma$ copies of $M$. 
Suppose adjacent molecules interact via a potential (energy) which is, by definition, 
a smooth function of two variables, say $f : M \times M \to \R$. Then the total 
energy of the crystal could be thought of as the (infinite !) sum of copies of $f$ 
over all adjacent pairs of sites $(\gamma, \gamma')$, called edges of $\Gamma$ 
(with six edges at each site)~:
$$
F(x= (x_\gamma)) = \sum_{{\rm edges}(\gamma, \gamma')}  f(x_\gamma, x_{\gamma'}).   
\eqno{(*)}
$$
Such an $F$ is clearly almost everywhere infinite for non-trivial $f$ but
its gradient (differential) is obviously well defined and finite at
all points $x$ in $X$. Thus one may speak of the critical points of $F$,
also called the stationary states of the crystal. Observe that the
set $S$ of stationary states make a closed subset in $X$ invariant under
the obvious (shift) action of $\Gamma$ on $X$. The basic question is that
of evaluating the entropy of the dynamical system $(S,\Gamma)$ in terms
of the topology of $M$ and/or some generic features of $f$. One still does
not have a satisfactory criterion for non-vanishing of this entropy
except for a few specific cases, such as the discretized geodesic flow
on a Riemannian manifold for instance, but one can give a lower bound
on the asymptotic distribution of the critical values of F as follows.\\

Exhaust $\Gamma$ with some standard subsets $\Omega_i$, e.g.~by concentric
cubes of edge size $2i$, and let $F_i$ denote the "restriction" of $F$ to $\Omega_i$, 
that is, the sum of the terms in $(*)$ corresponding to edges in $\Omega_i$. This sum
is regarded as a function on $M^{\Omega_i}$, call it $F_i : M^{\Omega_i} \to \R$. It
is further normalized by letting $F'_i = 1/\card(\Omega_i) \; F_i$. The functions
$F'_i$ take values in a fixed interval, namely in $[ f- = 6 \inf (f), f+ = 6 \sup (f)]$. 
We count the number $\#_i (I)$ of critical values in each subinterval $I$ of $[f-, f+]$, 
and set
$$
\cri_i(I) = \frac{1}{\card (\Omega_i)} \log \#_i(I).
$$
The purpose of this paper is to provide a Morse theoretic lower bound for 
$\liminf \cri_i(I)$, $i \to \infty$, in terms of a certain (strictly positive !) concave function 
(entropy) on $[f-,f+]$ capturing the homological behaviour of functions $F'_i$ for 
$i \to \infty$.\\

\noindent
{\bf Aknowledgements} We wish to thank the referee for having noticed some 
mistakes in the text. 

\section{Framework}\label{framework}

Consider a countable group $\Gamma$ endowed with a left-invariant
metric $d : \Gamma \times \Gamma \to \R^+$. Given a finite subset $\Omega 
\subset \Gamma$, its cardinality is denoted by $|\Omega|$. The set of finite 
subsets of $\Gamma$ is denoted by $B(\Gamma)$. The distance $d$ is 
extended to a map $d : B(\Gamma) \times B(\Gamma) \to \R^+$ (not a distance) 
as follows~: 
$$d(\Omega, \Omega') = \inf \{d(\gamma, \gamma') ; \gamma \in \Omega , 
\gamma' \in \Omega'\}.$$   
Given a nonnegative number $N$, the $N$-boundary and the $N$-interior 
of $\Omega \in B(\Gamma)$ are the sets  
$$\begin{array}{lcl}
\partial_N \Omega & = & \{ \gamma \in \Gamma ; d(\gamma, \Omega), 
d(\gamma , \Gamma - \Omega) \leq N\}, \\
{\rm int}_N\Omega & = & \Omega - \partial_N \Omega.
\end{array}$$ 
When reference to $N$ is clear the set ${\rm int}_{N/2}\Omega$ will be denoted 
by $\widetilde{\Omega}$. Given $\Omega_o$ and $\Omega$, two finite subsets of 
$\Gamma$, we denote their amenability ratio by $\alpha(\Omega, \Omega_o)$, that is
$$\alpha(\Omega, \Omega_o) = \frac{|\partial_{D_o} 
\Omega|}{|\Omega|},$$ where $D_o = \sup \{d(\gamma, \gamma') ; 
\gamma, \gamma' \in \Omega_o\}$ is the diameter of $\Omega_o$.
Let us recall that a countable group is said to be {\sl amenable} if 
it admits an amenable sequence $(\Omega_i)$, i.e.~an increasing 
sequence of finite subsets exhausting $\Gamma$ such that for any nonnegative 
number~$N$,
$$\lim_{i \to \infty} \frac{|\partial_N \Omega_i|}{|\Omega_i|} = 0.$$ \\

Let $X$ be a compact topological space endowed with a $\Gamma$-action 
$\rho :\Gamma \times X \to X : (\gamma , x) \mapsto \gamma x$.\\

Let $f_o : X \to \R$ be any continuous function with $f_o(X) = [0,1]$.
For $\Omega$ in $B(\Gamma)$, we define the average of $f_o$ 
along $\Omega$ to be the function
$$
f_\Omega (x) = \frac{1}{|\Omega|} \sum_{\gamma \in \Omega} 
f_o(\gamma^{-1}x).
$$

\section{Product-like actions}

We will impose on the group action $\rho$ a restrictive assumption of 
homological nature  expressing abundance of multiplicative structure. 
It will ensure that the homological measure defined in the next 
section will have a well-defined exponential growth. Its statement requires 
the introduction of some elements of notation.  \\ 

Let $F$ be a field and let $H^\star(X;F)$ denote the singular cohomology 
of $X$ with coefficients in $F$. Given a finite-dimensional subalgebra $A\subset    
H^\star(X;F)$ and a finite subset $\Omega$ of $\Gamma$, we denote by 
$A_{\Omega}$ the (finite-dimensional) subalgebra of $H^\star(X;F)$ generated 
by the translates of $A$ along $\Omega$, i.e.
$$A_\Omega = {\rm Alg} \; \Bigl< \; \bigoplus_{\gamma \in \Gamma} \gamma_* 
A \; \Bigr>,$$ where $\gamma_* = (\gamma^{-1})^*$ denotes the induced (left) action 
of $\gamma$ on $H^\star(X;F)$.   

\begin{ass}\label{ass}{\rm
There exists a nontrivial subalgebra $\a\subset H^\star(X;F)$ for which any finite-dimensional 
subalgebra $A \subset \a$ admits a number $N = N(A) \geq 0$ such that if 
$\Omega, \Omega' \in B(\Gamma)$ satisfy $d(\Omega, \Omega') > N(A)$, then 
the cup product map is injective 
\begin{equation}\tag{$\times$}
A_\Omega \otimes A_{\Omega'} \hookrightarrow A_{\Omega \cup \Omega'} : a \otimes a' 
\mapsto a \wedge a'.
\end{equation}
}\end{ass}

\begin{rmk}{\rm \aref{ass} The word {\em nontrivial} in the statement above should be given
the meaning that the algebra $\a$ contains some nonzero finite-dimensional  algebra.
}\end{rmk}
When this assumption is satisfied, the action $\rho$ is said to be a {\sl product-like 
action}, in reference to the following example. 

\begin{ex}\label{products}{\rm (Products)
Let $M$ be a manifold. Consider $X = M^{\Gamma}$, the infinite product of $\Gamma$ 
copies of $M$, or equivalently, the set of maps 
$$\Gamma \to M : \gamma \to x_\gamma,$$
with the topology of pointwise convergence (or product topology). 
\aref{ass} holds. Indeed, an algebra $\a$ satisfying the condition  
($\times$) is the direct limit of the direct system of subalgebras 
described hereafter. To each finite subset $\Omega \subset \Gamma$ is 
associated a finite-dimensional subalgebra $A(\Omega)$ of $H^\star( M^\Gamma; 
F)$~:    
$$A(\Omega) = p_\Omega^* ( H^\star(M^\Omega;F) ),$$
where $p_\Omega : M^\Gamma \to M^\Omega$ is the canonical projection. 
When $\Omega' \subset \Omega$, there is a map $p_{\Omega, \Omega'} : 
M^\Omega \to M^{\Omega'}$ and hence a pullback $i_{\Omega, \Omega'} : 
A(\Omega') \to A(\Omega)$. The algebra 
$$\a = \varinjlim_{\Omega \in B(\Gamma)} A(\Omega)$$ 
is a subalgebra of $H^\star( M^\Gamma; F)$ that satisfies ($\times$). Indeed, let $A 
\subset \a$ be a finite-dimensional subalgebra. There exists a finite subset 
$\Omega_o \subset \Gamma$ such that $A \subset A(\Omega_o)$. Given $\Omega$, 
the space $A_\Omega$ is contained in $A(\Omega \cdot \Omega_o)$ and, as 
the K\"unneth formula implies, it suffices to use $N(A) =  \diam (\Omega_o\cdot 
\Omega_o^{-1})$.   
}\end{ex}    
A class of examples of $\Gamma$-spaces not of the product type, but enjoying the 
product-like property is described below in \sref{NPE}. 

\begin{rmk}{\rm \aref{ass} is not satisfied when $X$ is a manifold, or 
when $H^\star(X;F)$ has finite rank.  
}\end{rmk}

\section{Homological measure of thickened level sets}

We will define homological invariants associated to a continuous function 
$f_o$ on $X$. They can be interpreted, roughly speaking, as a measure 
of the amount of cohomology supported in the various thickened level sets 
of the averages of $f_o$ over the finite subsets of $\Gamma$, and therefore 
could be called {\em homological measures} of slices. If $X$ was a manifold 
and if $f_o$ wasa smooth function, these invariants would provide a measure 
for the number of "homologically-detectable" critical points of the various 
functions $f_\Omega$ located in the various thickened level sets of 
$f_\Omega$ (cf.~\sref{Morse}). The real purpose is to consider the exponential 
growth of this invariant as the finite subset becomes large. The resulting object 
will depend upon two variables~: the level and the normalized degree in 
cohomology. It is called hereafter the {\sl homological entropy} of the function 
$f_o$, in analogy with the traditional entropy of an observable (\cite{L}).  
The entropy is well-defined provided these invariants satisfy certain properties. 
Classically these properties are submultiplicativity and $\Gamma$-invariance. 
In contrast, the homological measure is invariant as well (\lref{inv}), but 
{\bf super}multiplicative (\lref{superm}). To define entropy in this situation 
necessitates the introduction of an additional assumption on the group, called 
here {\sl tileability} (cf.~\sref{OW}).    

\begin{ntt}
If $a$ is a cohomology class and if $O$ is an open subset of $X$, the 
expression $\supp a \subset O$ ("$a$ is supported in $O$"), means 
that for some open set $O'$ such that $X = O \cup O'$, the restriction of 
$a$ to $O'$ vanishes. Observe that if $\supp a \subset O$ and $\supp b 
\subset U$ then $\supp a \wedge b \subset O \cap U$. 
\end{ntt}

Given $A \subset \a$ a finite-dimensional subalgebra, an open set $O$, 
a non-negative number $\ell$, a positive number $\nu$ and a finite subset 
$\Omega \subset \Gamma$, we consider the subalgebras      
$$\begin{array}{rcl}
\Ho{\star}{O} & = & \{ a \in A_{\widetilde \Omega} ; \supp a \subset 
O\}, \\
\Ho{\ell, \nu}{O} & = & \{ a \in A_{\widetilde \Omega} ; \supp a \subset 
O \; \& \; (\ell - \nu) |\Omega| < \deg a < (\ell + \nu) |\Omega| 
\}.
\end{array}$$
Recall the convention that $\widetilde{\Omega}$ denotes $\Int_{N/2} \Omega$,
where $N = N(A)$ is the number associated to $A \subset \a$ from \aref{ass}. 
The inequalities involving the degree of $a$ have to be verified by each component 
of pure degree. Finally, the open sets considered hereafter will be sublevel, superlevel or 
thickened level sets of the function $f_\Omega$, typically~:
$$
O = f_{\Omega}^{-1} (- \infty , c + \delta ) \text{ or } 
f_{\Omega}^{-1} (c - \delta,  +\infty ) \text{ or } f_{\Omega}^{-1} 
(c - \delta , c + \delta ), 
$$
for some $c \in [0,1]$ and $\delta > 0$. Then consider the map

$$\begin{array}{rcl}
\ph{\ell, \nu} & : & \Hm{\ell, \nu} \quad \rightarrow  \\
&&\\
& & \Hom \biggl( \Hp{\star} , \Hpm{\star} \biggr)\\
\end{array}$$
$$\Bigl(\ph{\ell, \nu}(a)\Bigr)(b) =  a \wedge b.$$
Its rank is denoted hereafter by $$\bi{\ell, \nu} = \rk \ph{\ell, 
\nu}.$$ 

\begin{dfn} $\bi{\ell, \nu} $ is called the $(\ell - \nu, \ell + \nu)$-th homological 
measure of the thickened level $f_\Omega^{-1}(c - \delta, c + \delta)$ 
with respect to $A_\Omega$.
\end{dfn}

\section{Properties of the homological measure}

We prove in this section the two properties -- supermultiplicativity 
and $\Gamma$-invariance -- necessary to obtain a well-defined 
homological entropy.

\begin{lem}\label{superm} The map $B(\Gamma) \to \R : \Omega \mapsto \bi{\ell, 
\nu}$ is supermultiplicative. More generally, let $\Omega, \Omega' \in  
B(\Gamma)$ be disjoint finite subsets, let $c, c' \in [0,1]$ and let $\ell,\ell'\geq 
0$, then  
\begin{equation}
\biun{\alpha \ell + (1 - \alpha) \ell', \nu} \geq \bi{\ell, \nu} \cdot
\bipr{\ell', \nu},
\end{equation}
where $\alpha = \frac{|\Omega|}{|\Omega \cup \Omega'|}$ and thus $1 - 
\alpha = \frac{|\Omega'|}{|\Omega \cup \Omega'|}$.
\end{lem}  
  
\Pf The argument relies on the few simple observations listed below~:

\begin{enumerate}

\item[-] Let $\Omega$ and $\Omega'$ be disjoint finite subsets 
of $\Gamma$, then
$$f_{\Omega \cup \Omega'} = \tfrac{|\Omega|}{|\Omega \cup \Omega'|} 
f_\Omega + \tfrac{|\Omega'|}{|\Omega \cup \Omega'|} f_{\Omega'} = 
\alpha f_\Omega + (1 - \alpha) f_{\Omega'}.$$
Thus $f_{\Omega \cup \Omega'}^{-1}(\alpha I + (1 - \alpha) I') \supset 
f_{\Omega}^{-1}(I) \cap f_{\Omega'}^{-1}(I')$ for intervals $I$ 
and $I'$. 

\item[-] If a class $a$ is supported in $f_{\Omega}^{-1}(I)$ and if a class 
$a'$ is supported in $f_{\Omega'}^{-1}(I')$, then the class $a \wedge a'$ 
is supported in $f_{\Omega}^{-1}(I) \cap f_{\Omega'}^{-1}(I') \subset 
f_{\Omega \cup \Omega'}^{-1}(\alpha I + (1 - \alpha) I')$.

\item[-] If $\Omega$ and $\Omega'$ are disjoint then the distance between 
$\widetilde{\Omega}$ and $\widetilde{\Omega}'$ is greater than  $N$ and 
$\widetilde{\Omega} \cup \widetilde{\Omega}' \subset \widetilde{ \Omega 
\cup \Omega'}$. Therefore \aref{ass} provides us with an 
injective map 
$$A_{\widetilde{\Omega}} \otimes A_{\widetilde{\Omega}'} \to 
A_{\widetilde{\Omega \cup \Omega'}}.$$
This explains the choice of $\widetilde{\Omega}$ instead of $\Omega$.

\item[-] The degree of $a \wedge a'$ if the sum of the degree of $a$ 
and that of $a'$. Thus 
$$\left. \begin{array}{rcl} 
a \in H^{\ell, \nu}_{A_{\widetilde{\Omega}}} \\ 
a' \in H^{\ell', \nu}_{A_{\widetilde{\Omega}'}} 
\end{array} \right\} \quad \Rightarrow \quad 
a \wedge a' \in H^{\alpha \ell + (1 - \alpha) \ell', \nu}_{A_{\widetilde{\Omega 
\cup \Omega'}}}.$$  
\end{enumerate}
Combining the previous observations we obtain an injection~:

$$\Psi_{I, I'} : \Ho{\ell, \nu}{f_{\Omega}^{-1}( 
I)} \otimes \Hopr{\ell', \nu}{f_{\Omega'}^{-1} (I')} \to \Houn{\alpha \ell + 
(1 - \alpha) \ell', \nu}{f_{\Omega \cup \Omega'}^{-1} (\alpha I + (1 - 
\alpha) I')}.$$ 
Now consider the following sequence of maps. We will abbreviate $\alpha 
c + (1 - \alpha) c'$ to $\tilde{c}$ and $\alpha \ell + (1 - \alpha) \ell'$ to 
$\tilde{\ell}$. 

$$\begin{array}{c}
\Hm{\ell, \nu}\bigotimes\Hmpr{\ell', \nu} \\
\downarrow\\
\Hom\biggl( \Hp{\star},\Hpm{\star} \biggr) \bigotimes \\
\Hom \biggl( \Hppr{\star},\Hpmpr{\star} \biggr) \\
\downarrow \\
\Hom \biggl( \Hp{\star}\otimes\Hppr{\star}, \\
\Hpm{\star}\otimes\Hpmpr{\star} \biggr) \\
\downarrow \\
\Hom\biggl( \Hpun{\star},\Hpmun{\star} \biggr).
\end{array}$$
The first arrow stands for the map $\ph{\ell, \nu} \otimes \phpr{\ell', 
\nu}$. The second arrow is a classical isomorphism, indeed, $$\Hom(A, B) 
\otimes \Hom(C, D) \simeq \Hom(A \otimes C, B \otimes D)$$ for 
finite-dimensional vector spaces $A,B,C,D$. The third one 
is the injection induced by a choice of complementary subspaces to the 
images of the maps $\Psi_{(c - \delta, + \infty), (c' - \delta, + \infty)}$ 
and $\Psi_{(c - \delta, c + \delta), (c' - \delta, c' + 
\delta)}$ respectively. We will denote the composition of second and third 
map by $\Phi$. There is also another sequence obtained from composing 
$$\Psi = \Psi_{(- \infty, c + \delta), (- \infty, c' + \delta)}$$ with 
$$\phun{\alpha \ell + (1 - \alpha) \ell', \nu}.$$ These two sequences commute~:  

$$\Phi \circ \bigl[\ph{\ell, \nu} \otimes \phpr{\ell', \nu} \bigr] = 
\phun{\alpha \ell + (1 - \alpha) \ell', \nu} \circ \Psi.$$
Since $\Phi$ and $\Psi$ are both injective, this implies that  
$$ \rk \Bigl( \ph{\ell, \nu} \otimes \phpr{\ell', \nu} \Bigr) \leq \rk 
\phun{\alpha \ell + (1 - \alpha) \ell', \nu}$$

\cqfd

\begin{lem}\label{inv} The map $$B(\Gamma) \to \R : \Omega \to \bi{\ell, \nu}$$
is $\Gamma$-invariant.
\end{lem}

\Pf
The proof follows from the simple facts stated below. Let $\Omega \in 
B(\Gamma)$, let $\gamma_o \in \Gamma$ and let $I \subset [0,1]$ be any 
interval. Then   

\begin{enumerate}
\item[-] $f_{\gamma_o \Omega}^{-1} (I) = \gamma_o f_\Omega^{-1} (I)$.
\item[-] $A_{\gamma_o \Omega} = (\gamma_o)_* A_\Omega$.
\item[-] $\widetilde{\gamma_o \Omega} = \gamma_o \widetilde{\Omega}$. 
\item[-] $(\gamma_o)_*$ induces an isomorphism between 
$\Ho{\ell, \nu}{f_\Omega^{-1} (I)}$ and $H_{{\gamma_o}_* 
A_{\widetilde{\Omega}}}^{\ell, \nu}(\gamma_o f_{\Omega}^{-1} (I))$. 
\end{enumerate}

\cqfd

\section{Superadditive Ornstein-Weiss Lemma for ti\-lea\-ble groups}\label{OW}

The $\ell$-th Betti number entropy of $f_o$ will be defined from the exponential 
growth, with respect to the index $i$, of the sequence $(\bil{\ell, \nu}{i})_{i \geq 1}$, 
where $\Omega_i$ is an amenable sequence in $\Gamma$, that is to say, from 
the limit $$\lim_{i \to \infty} \frac{\ln ( \bil{\ell, \nu}{i} )}{|\Omega_i|},$$ 
when it exists. \lref{superm} implies that the map $\Omega \mapsto \ln ( 
\bi{\ell, \nu} )$ is superadditive on disjoint sets, while the Ornstein-Weiss  
lemma \cite{O-W} provides convergence of the sequence of averages 
$h( \Omega_i ) / |\Omega_i|$ under the hypotheses that the map $h : 
B(\Gamma) \to \R^+$ is $\Gamma$-invariante and {\bf sub}additive. The proof 
of the Ornstein-Weiss Lemma requires the construction of 
$\varepsilon$-quasi-tilings that any amenable group admits. In contrast, a 
proof of the superadditive version of this lemma seems to necessitate the 
construction of {\sl disjoint} such tilings which might not exist in general (although we 
do not know of any counterexample). Whence the following definition.     

\begin{dfn}\label{tileable} An amenable group $\Gamma$ is said to be tileable 
if it admits a tiling amenable sequence, that is, an amenable sequence $(\Omega_i)$ 
such that given $\varepsilon > 0$ and a subsequence $(\Omega_{i_n})$ of $(\Omega_i)$, 
there exists a finite subsequence of $(\Omega_{i_n})$, denoted $\Omega_1, \ldots, \Omega_s$, 
such that any finite subset $\Omega$ with sufficiently large amenablility ratios $\alpha(\Omega,
\Omega_j), j = 1, \ldots, s$ can be disjointly $\varepsilon$-tiled by translates of the $\Omega_j$'s,
i.e.~there exists center $\gamma_{j,k}, 1 \leq j \leq s, 1\leq k \leq r_j$ in $\Gamma$ such that  
\begin{enumerate}
\item[-] $\gamma_{j,k} \Omega_j \subset \Omega$,
\item[-] $\gamma_{j,k} \Omega_j \cap \gamma_{j',k'} \Omega_{j'} = 
\emptyset$ for $(j,k) \neq (j',k')$,
\item[-] $|\cup_{j,k} \gamma_{j,k} \Omega_{j}| \geq (1 - \varepsilon) |\Omega|$.
\end{enumerate}
\end{dfn}

\begin{exs}{\rm Weiss introduces in \cite{W} the notion of {\em monotileable amenable 
groups}. Those are groups admitting an amenable sequence $(\Omega_i)$ consisting 
of monotiles. This means that for each index $i$ there exists a set $C_i \subset \Gamma$ 
for which the various translates $\Omega_i c$ of $\Omega_i$ along $C_i$ form a partition 
of $\Gamma$.
Such groups belong to the class of tileable amenable groups. Moreover, Weiss proves 
that any residually finite amenable group is monotileable, implying that the following 
amenable groups are also tileable.
\begin{enumerate}
\item[-] Abelian and solvable groups.
\item[-] Amenable linear groups, i.e.~linear groups not containing $F_2$ as a subgroup.
\item[-] Grigorchuk's groups of intermediate growth.\footnote{Grigorchuk, Rotislav I.,
Degrees of growth of finitely generated groups and the theory of invariant means. 
{\em Izv. Akad. Nauk SSSR Ser. Mat.} {\bf 48} (1984), no. 5, 939--985.}
\end{enumerate}
}\end{exs}

\begin{lem}\label{OWL}(Superadditive Ornstein-Weiss lemma) Let $\Gamma$ be a 
tileable amenable group. Let $h$ be a nonnegative function defined on $B(\Gamma)$ 
and satisfying the following two conditions~: 
\begin{enumerate}
\item[-] superadditivity~: $h(\Omega \cup \Omega') \geq h(\Omega) + 
h(\Omega')$ for disjoint subsets $\Omega$ and $\Omega'$,
\item[-] $\Gamma$-invariance~: $h(\gamma \Omega) = h(\Omega)$ 
for any $\gamma \in \Gamma$.
\end{enumerate}
Then, given a tiling amenable sequence $(\Omega_i)$, the following limit 
exists
$$\lim_{i \to \infty} \frac{h(\Omega_i)}{|\Omega_i|}.$$
\end{lem}

\begin{rmk}{\rm Observe that under the hypotheses of the previous lemma, 
the limit is independent of the choice of a tiling amenable sequence in $\Gamma$.
}\end{rmk}

\Pf Let $\varepsilon > 0$ and let $(\Omega_i)$ be a tiling amenable sequence. Extract 
a subsequence $\Omega_{i_1}, \ldots, \Omega_{i_s}$ with which we can 
$\varepsilon$-tile any element $\Omega_i$ of the initial sequence with sufficiently 
large index. Suppose also that if $h^+$ stands for the limsup of the sequence
$h(\Omega_i)/ |\Omega_i|$, then $h(\Omega_{i_j})/ |\Omega_{i_j}| \geq h^+ - 
\varepsilon$ for all $j$. Let $\gamma_{j,k}, 1 \leq j \leq s, 1\leq k \leq r_j$ 
denote the centers of a disjoint $\varepsilon$-tiling of $\Omega_i$ by translates of 
the $\Omega_{i_j}$'s and let $\Omega'_i = \cup_{j,k} \; \gamma_{j,k} \Omega_{i_j}$. 
Then    
$$\begin{array}{rcl}
\dfrac{1}{|\Omega_i|} h(\Omega_i) & \geq & \dfrac{1}{|\Omega_i|}  
\Bigr( h( \Omega'_i ) + h( \Omega_i - \Omega'_i) \Bigr) \geq 
\dfrac{1}{|\Omega_i|} h( \Omega'_i ) \\    
& \geq & \dfrac{1}{|\Omega_i|} \sum_{j,k} h(\Omega_{i_j}) \geq 
\dfrac{1}{|\Omega_i|} \sum_{j,k} (h^+ - \varepsilon) |\Omega_{i_j}| \\  
& \geq & \dfrac{1}{|\Omega_i|} (h^+ - \varepsilon) (1 - \varepsilon) 
|\Omega_i| = (h^+ - \varepsilon) (1 - \varepsilon).
\end{array}$$
Hence $$\liminf_{i \to \infty} \dfrac{1}{|\Omega_i|} h(\Omega_i) \geq 
(h^+ - \varepsilon) (1 - \varepsilon).$$ Since this holds for arbitrary 
$\varepsilon$, the limit $\lim_{i \to \infty} h(\Omega_i) / |\Omega_i|$ exists.

\cqfd

\section{Homological entropy of functions}
        
Let $(\Omega_i)_{i\geq 1}$ be a tiling amenable sequence in the tileable 
amenable group $\Gamma$ and consider the exponential growth of the sequence 
$\bil{\ell, \nu}{i}$~: 
\begin{equation}\label{limit}
b^{\ell, \nu}_{A}(c,\delta) = \lim_{i \to \infty} \dfrac{1}{|\Omega_i|} 
\ln \Bigl[\bil{\ell, \nu}{i}\Bigr].  
\end{equation}
As implied by \lref{superm}, \lref{inv} and \lref{OWL}, this limit indeed exists.
Observing that  the function $b^{\ell, \nu}_{A}(c,\delta)$ is increasing in $\nu$ and 
$\delta$, we let $\delta$ and $\nu$ approach~$0$~:

\begin{equation*}
b^{\ell}_{A}(c) = \lim_{\nu \to 0} \lim_{\delta \to 0} b^{\ell, 
\nu}_{A}(c, \delta).  
\end{equation*}
Independence on $A$ is obtained by considering the supremum over all 
possible choices of a finite-dimensional subalgebra $A$ of $\a$~:
\begin{equation*}
b^\ell(c) = \sup\{ b^\ell_{A}(c) ; A \subset \a \;\; \& \; \Dim A < \infty\}. 
\end{equation*}
 
\begin{dfn}
The function $b^\ell : [0,1] \to \R : c \mapsto b^\ell(c)$ is called the 
$\ell$-th Betti number entropy of $f_o$.
\end{dfn}    

\begin{rmk}{\rm One may also define the {\em sum of the Betti number entropy of 
$f_o$} by the same process except that the cohomological degree is not restricted. The 
two functions are related as follows~:
$$b(c) = \sup_\ell b^\ell(c).$$ 
(This is a consequence of the general fact that the exponential growth of the sum of two 
sequences coincides with the exponential growth of the maximum sequence.)
}\end{rmk}

\begin{rmk}{\rm The condition that $\Gamma$ be tileable is only used to guarantee
existence of the limit (\ref{limit}).
}\end{rmk}

\section{Relation with classical Morse theory}\label{Morse}

This section is devoted to showing how, in the setting of a manifold $M$ 
endowed with a Morse function $f:M \to \R$, the sum of the Betti number entropy 
essentially coincides with the sum of the Betti numbers of $M$ and provides a lower 
bound for the number of critical points of $f$ (cf.~\pref{MI}). \\

Let $M$ be a connected, closed, oriented smooth manifold. Consider a
Morse function $f$ on $M$. Let $F$ be a field. We recall that if $O$ is some 
subset of $M$ and if $a$ is a cohomology class in $H^\star(M;F)$, the expression 
$\supp a \subset O$ means that $a|_{O'} = 0$ for some open subset $O'$ 
containing $M - \Int O$. We denote by $\Crit_c(f)$ the set of critical points of $f$ at 
level $c$. Define

$$\begin{array}{l} 
H^\star(O) = \Bigr\{ a \in H^\star(M;F) ; \supp a \subset O \Bigl\}, \\
\end{array}$$
$$\begin{array}{rll} 
\varphi_{c,\delta} : H^\star(f^{-1}(- \infty, c + \delta)) & \to & \Hom \Bigl( H^\star(f^{-1}
(c - \delta, + \infty)), H^\star(f^{-1}(c - \delta, c + \delta))\Bigr), \\
a & \mapsto & \Bigl[ b \mapsto a \wedge b\Bigr],
\end{array}$$
$$\begin{array}{l} 
b(c, \delta) = \rk \varphi_{c,\delta}, \\
\end{array}$$
$$\begin{array}{l} 
b(c) = \lim_{\delta \to 0} b(c, \delta).
\end{array}$$
We might sometimes denote $b(c, \delta)$ by $b(c - \delta, c + \delta)$ or consider 
$b(I)$ when $I$ is some interval.

\begin{prop}\label{MI}\hspace{1cm}
\begin{enumerate}
\item[(a)] $\displaystyle{\sum_{c \in \R}} b(c) = SB(M)$.
\item[(b)] $\displaystyle{b(c) \leq \Crit_c(f)}$.
\end{enumerate}
\end{prop}

\Pf \\
(a) The main ingredient is the specific version of Poincar\'e duality mentioned 
below in \sref{PDcom}. Indeed, if $a \in H^\star(M; F)$,  define
$$c_a = \inf \{c ; \supp a \subset f^{-1}(- \infty, c)\}.$$
Since for all $\delta > 0$ the restriction of $a$ to $f^{-1}(c_a - \delta, + \infty)$ 
does not vanish, there exists a class $b$ with $\supp b \subset f^{-1}(c_a - \delta, 
+ \infty)$ such that $a \wedge b \neq 0$ (cf.~\pref{PDcom-prop}). Hence $a$ 
provides a contribution to $b(c_a)$. Here follows a more precise argument taking 
into account the following difficulty~: there might exist two classes that have same 
$c_a$ and are independent, but who do not generate a $2$-dimensional space of 
classes with same $c_a$. \\

\noindent
Decompose the range $I$ of $f_o$ into intervals as follows~:
$$\displaystyle{I \subset \bigcup_{k = 1}^K (I_k = [a_k, a_{k+1}))} \qquad a_k < 
a_{k+1}.$$
If $J \subset \R$, let $r_J : H^\star(M;F) \to H^\star(f^{-1}(J); F)$ denote the restriction 
map. Then consider the following increasing sequence of subspaces
$$\{0\} = \Ker r_{I_1 \cup ... \cup I_K}\subset ... \subset \Ker r_{I_{K-1} \cup I_K } 
\subset \Ker r_{I_K} \subset H^\star(M;F).$$
Choose a corresponding sequence of spaces $V_1, ... , V_K$ such that
$$\Ker r_{I_k \cup ... \cup I_K} \oplus V_k = \Ker r_{I_{k+1} \cup ... \cup I_K} 
\qquad k= 1, ...,K.$$
(For $k=K$, we mean $\Ker r _{I_K} \oplus V_K = H^\star(M;F)$.)
If $0 \neq a \in V_k$ then $c_a \in I_k$. Hence $b(I_k) = \rk V_k$ and 
\begin{equation}\label{M<}
\displaystyle{\sum_{k=1}^K b(I_k) = SB(M)}.
\end{equation}
This is true for arbitrarily fine subdivisions of $I$. Now observe that
for each $c$, either $b(c) = 0$, in which case $b(c - \delta, c + \delta) = 0$ for 
all sufficiently small $\delta > 0$, or $b(c - \delta, c + \delta) \neq 0$ for all $\delta > 0$.
Relation (\ref{M<}) implies that there are finitely many numbers $c$ with $b(c) \neq 0$. 
So  $\sum_kb(I_k)$ is constant, equal to $\sum_c b(c)$, for all sufficiently fine 
subdivisions of $I$. \\

\noindent
(b) Given $a \in H^\star(M;F)$, $c_a$ must be a critical value, otherwise we 
would be able to move $f^{-1}(- \infty, c + \delta)$ below level $c - \delta$
by an ambient isotopy, disjointifying the supports of the classes $a$ and $b$.
In consequence, $a \wedge b$ would vanish. This alone implies (b) when $b(c) 
\leq 1$ for all $c$.
We will argue that if $b(c) = 2$ then $f$ cannot have a single critical point 
at level $c$ (the general case can be handled in a similar way). Suppose on the 
contrary that $\{x\} = \Crit_c(f)$. Let $x_1, ..., x_m$ be coordinates on $M$, centered at $x$, 
for which $f$ has the canonical form 
$$f(x) = - x_1^2 - ... - x_n^2 + x_{n+1}^2 + ... + x_m^2.$$ 
Let $a_1$ and $a_2$ be two independent classes with $c_{a _1}= c_{a_2} = c$. 
Consider piecewise smooth cycles $\alpha_1$ and $\alpha_2$ representing 
their Poincar\'e dual homology classes. We will make the following assumptions 
on $\alpha_i$, $i = 1, 2$~:

\begin{enumerate}
\item[-] $\alpha_i$ is supported in $f^{-1} (-\infty, c]$,
\item[-] $\alpha_i \cap f^{-1}(c) = \{x\}$,
\item[-] $x$ is a regular value of (each of the simplices composing) $\alpha_i$,
\item[-] $\alpha_i$ intersects the local unstable manifold ${\mathcal W}^u(x) 
= \{x ; x_{1} = ... = x_n = 0\}$ of $x$ transversely at $x$. 
\end{enumerate}
It is long but not difficult to verify that these hypotheses are not restrictive. \\

Now, we will show that the degree of $\alpha_i$ must equal the index $n$ of $x$.
The degree of $\alpha_i$ can certainly not exceed $n$, otherwise $\alpha_i$ 
would not be supported in $f^{-1} (-\infty, c]$. If the degree of $\alpha_i$ was less than $n$, 
one could slide $\alpha_i$ down the stable manifold of $x$ (in a direction transverse to 
that of $T_x\alpha_i$) below level $c$. \\

Then, one can subdivide all the simplices of $\alpha_i$ containing $x$ in such 
a way that $x$ lies in the interior of each simplex to which it belongs and that each such
simplex can be isotoped to a fixed simplex that coincides with the stable manifold 
${\mathcal W}^s(x)$ of $x$ in a neighborhood of $x$. Thus, 
$$\alpha_i = f_i^0 \sigma_i^0 + \sum_{j \geq 1} f_i^j \sigma_i^j,$$ 
where $f_i^0, f_i^j \in F$, where $\sigma_i^0$ is a piece of ${\mathcal W}^s(x)$ 
containing $x$ and where $\sigma_i^j$ avoids $x$. It follows that $f_1^0 \alpha_2 - 
f_2^0 \alpha_1$ vanishes near $x$, hence that $c_{f_1^0 a_2 - 
f_2^0 a_1} < c$. So $a_1$ and $a_2$ do not generate a space contributing to $b(c)$, 
a contradiction.

\cqfd

\begin{rmk}{\rm The previous lemma implies the Morse theoretic lower bound announced 
in the introduction. Referring to the notation used therein, one observes that the previous 
proof implies in particular that if $I \subset \R$ is some interval, then $\Crit_I(F'_i) \geq 
b_{A, \Omega_i}(I)$, where $A \simeq H^*(M; F)$. Hence
$$\liminf_{i\to \infty} \cri_i(I) \geq b_A(I).$$
(The function $F'_i$ defined in the introduction does not quite coincide with the function
$f_{\Omega_i}$ defined in \sref{framework}, but the difference will not affect the asymptotic 
behavior of the objects considered here). Moreover, as proved later on, the function $b_A(I)$ 
is concave and strictly positive.
}\end{rmk}

\section{Concavity of the entropy}

The above-defined function $b : \R^+ \times [0,1] \to \R$ is concave. This follows mainly from
\lref{superm}, with a slight help from the following fact. 

\begin{lem}\label{usc} The function $b$ is upper semi-continuous.
\end{lem}

\Pf Let $(\ell_k, c_k)$ be a sequence converging to some pair $(\ell, c)$ in $\R^+ \times [0,1]$. 
Since $b^{\ell, \nu}_{A, \Omega} (c, \delta)$ is increasing \wrt both intervals $(\ell - \nu, \ell + \nu)$  
and $(c - \delta, c + \delta)$, 
$$b^{\ell, \nu}_{A, \Omega} (c, \delta) \geq b^{\ell_k, 
\frac{\nu}{2}}_{A, \Omega} (c_k, \tfrac{\delta}{2})$$ 
for sufficiently large $k$ and for all $A$ and $\Omega$. Hence $b^\ell (c) \geq 
b^{\ell_k} (c_k)$. Thus $b^\ell (c) \geq \limsup_{k \to \infty} b^{\ell_k} (c_k)$.

\cqfd
    
\begin{prop} The function $b$ is concave. That is to say, for any $\ell, \ell'\in 
\R^+$, any $c, c' \in [0,1]$, and any $\alpha \in [0,1]$, 
\begin{equation}\label{cvxty}
b^{\alpha \ell + (1 - \alpha) \ell'} (\alpha c + (1 - \alpha) c') 
\geq \alpha \; b^{\ell} (c) + (1 - \alpha) \; b^{\ell'} (c').
\end{equation}
\end{prop}

\Pf Let $(\Omega_i)$ be an amenable sequence. For each $i$, 
let $\Omega'_i = \gamma_i \cdot \Omega_i$ be disjoint from $\Omega_i$. 
Then the sequences $(\Omega'_i)$ and $(\Omega_i \cup \Omega'_i)$ are amenable as well. Besides, \lref{superm} 
implies that 
$$b^{\alpha \ell + (1 - \alpha) \ell', \nu}_{A, \Omega_i \cup 
\Omega'_i}(\alpha c + (1 - \alpha) c', \delta) \geq b^{\ell, \nu}_{A, 
\Omega_i} (c, \delta) \cdot b^{\ell', \nu}_{A, \Omega'_i} (c', \delta),$$
with $\alpha = \frac{1}{2}$. Hence
\begin{equation}\label{dyadic-cvxty}
b^{\tfrac{1}{2} \ell + \tfrac{1}{2} \ell'} \Bigl(\tfrac{1}{2} c + 
\tfrac{1}{2} c' \Bigr) \geq \tfrac{1}{2} b^{\ell} (c) + \tfrac{1}{2} 
b^{\ell'}(c'). 
\end{equation}  
This implies that the relation (\ref{cvxty}) holds for any dyadic rationnal $\alpha$. The result 
for arbitrary $\alpha$ follows from the upper semi-continuity of $b$ (\lref{usc}).  
\cqfd

\section{Nontriviality of the entropy for products}

Let $F$ be a field. Let $M$ be a closed $F$-orientable manifold, that is to say 
$H^m(M;F) \simeq F$ for $m = \Dim M$. In other words, either $M$ is 
orientable, or $F = {\mathbb Z}_2$. Let $f_o : X = M^\Gamma \to \R$ be a 
continuous function with range $[0,1]$. 

\begin{prop}\label{pos} The associated homological entropy of $f_o$ achieves a 
strictly positive value.  
\end{prop}
This results holds because a products $M^\Gamma$ inherits some Poincar\'e 
duality (cf. \lref{PD}) from the manifold $M$. 

\subsection{Poincar\'e duality on a closed orientable manifold}\label{PDcom}

Here follows the specific version of Poincar\'e duality that is needed 
below. 

\begin{prop}\label{PDcom-prop}
If $a$ is a class in $H^\star(M;F)$ whose restriction to the open set $O$ does not vanish, 
then there exists a class $b$ with support in $O$ such that $a \wedge b \neq 0$.  
\end{prop}

\Pf Let $a \in H^i(M; F)$ with $a|_O \neq 0$. Then there exists a homology 
class $\beta \in H_i(O; F)$ such that $<a, \beta> \neq 0$. Let $b$ be 
the Poincar\'e dual of $\beta$. Then $a \wedge b$ does not vanish 
since its evaluation on the fundamental class of $M$ coincides with $<a, 
\beta>$. Moreover, if $\beta$ is represented by a chain $c$, the class 
$b$ can be represented by a form whose support is contained in any given 
neighborhood of the image of~$c$. 
\cqfd

\subsection{Poincar\'e duality in a product $M^\Gamma$}

Let $id \in \Omega_o  \subset \Gamma$ be a finite subset and let $A = A(\Omega_o) 
= p_{\Omega_o}^*( H^\star(M^{\Omega_o}; F))$, where $p_{\Omega_o} : M^\Gamma 
\to M^{\Omega_o}$ is the canonical projection. Let also $N = N(A)$ (cf.~\aref{ass} and 
\eref{products}). If $\Omega \subset \Gamma$ is another finite subset, we can define the 
positive number    
$$
\delta_\Omega = \delta_\Omega (f_o, \Omega_o) = 4\; \sup \{ \bigl| 
f_\Omega (x) - f_\Omega (y) \bigr| ; x_\gamma = y_\gamma \text{ for } 
\gamma \in \Omega \cdot \Omega_o \}.
$$
Observe that $\delta_\Omega$ is decreasing in $\Omega$. In fact 
$\delta_\Omega$ approches $0$ as $\Omega$ becomes large (cf.~\lref{delta}).

\begin{lem}\label{PD}(Poincar\'e duality in $M^\Gamma$) Let $a$ belong 
to $A_\Omega$. Then there exists a level $c_a^{\Omega}$ and an element 
$b$ in $A_\Omega$ such that  

\begin{enumerate}
\item[-] $\supp a \subset f_\Omega^{-1}(- \infty , c_a^{\Omega} + 
\delta_\Omega)$, 
\item[-] $\supp b \subset f_\Omega^{-1}(c_a^{\Omega} - \delta_\Omega, + 
\infty)$,   
\item[-] $a \wedge b \neq 0$.
\end{enumerate}

\end{lem}

\Pf
The level $c_a^{\Omega}$ defined below obviously satisfies the first 
condition. 
\begin{equation}\label{level}
c_a^{\Omega} = \inf \{ c \in [0,1] ; \supp a \subset f_\Omega^{-1}(- \infty, 
c) \}.
\end{equation}
As explained below, existence of  the class $b$ follows from Poincar\'e duality 
in any finite product $M^\Omega$. Choose a point $o$ in $M^\Gamma$ and 
define new functions $g_\Omega : M^\Gamma \to [0,1]$  by $g_\Omega (x) = 
f_{\Omega} (\hat{x})$, with 
$$\hat{x}_\gamma = \Biggl\{ \begin{array}{ccc} x_\gamma & \text{ if } 
\gamma \in \Omega \cdot \Omega_o \\
o_\gamma & \text{ otherwise. }
\end{array}$$  
By definition of $\delta_\Omega$, 
$$\sup_{x \in M^{\Gamma}} \bigl| f_\Omega(x) - g_\Omega(x) \bigr| \leq 
\frac{\delta_\Omega}{4}.$$

Now observe that the restriction of $a$ to $g_\Omega^{-1} (c_a^\Omega - 
\frac{3}{4} \delta_\Omega, + \infty)$ does not vanish. Indeed, if it did 
vanish then $\supp a$ would be contained in $g_\Omega^{-1} (- 
\infty, c_a^\Omega - \frac{1}{2} \delta_\Omega)$ which itself is 
contained in $f_\Omega^{-1} (- \infty, c_a^\Omega - \frac{1}{4} 
\delta_\Omega)$. This contradicts the definition of~$c_a^\Omega$.\\

Since $g_\Omega$ depends only on the variables indexed by 
$\Omega \cdot \Omega_o$, the open set $g_\Omega^{-1} (c_a^\Omega - 
\frac{3}{4} \delta_\Omega, + \infty)$ coincides with the pullback 
$p_{\Omega \cdot \Omega_o}^{-1} (O)$ of some open subset $O$ of 
$M^{\Omega \cdot \Omega_o}$. Combined with the fact that $a = 
p_{\Omega \cdot \Omega_o}^* \ol{a}$ for some class $\ol{a}$ in 
$H^\star(M^{\Omega \cdot \Omega_o}; F)$, this implies that the restriction 
of the class $\ol{a}$ to $O$ does not vanish. Poincar\'e duality in 
closed orientable manifolds yields a class $\ol{b} \in H^\star( M^{\Omega 
\cdot \Omega_o}; F)$ with $\supp \ol{b} \subset O$ and such that $\ol{a} \wedge 
\ol{b} \neq 0$. The class $b = p_{\Omega \cdot \Omega_o}^{*} 
(\ol{b})$ satisfies the required conditions.

\cqfd

The following result implies that in \lref{PD} one can replace $\delta_\Omega$
by any given $\delta>0$ at the cost of considering only "large" $\Omega$'s.

\begin{lem}\label{delta} Let $(\Omega_i)_{i\geq 1}$ be an amenable sequence. 
Then the sequence $\delta_{\Omega_i}$ converges to $0$.
\end{lem}

\Pf Let $\delta > 0$. By (uniform) continuity of $f_o$, there exists a 
$\eta = \eta(\delta)> 0$ such that $\hat{d}(x,y) < \eta \Rightarrow |f_o(x) - 
f_o(y)| < \delta$. The symbol $\hat{d}$ denotes one of the following (compatible) 
metrics on $M^\Gamma$~: $$\hat{d}(x,y) = \sum_{\gamma \in \Gamma} 
\frac{d_o(x_\gamma, y_\gamma)}{\lambda^{|\gamma|}},$$ where $d_o$ is 
some Riemannian metric on $M$, where $\lambda$ is some fixed number in 
$(1, +\infty)$ and where $|\gamma| = d(id, \gamma)$. In particular $\hat{d}(x,y) < 
\eta$ when sufficiently many components of $x$ and $y$ coincide. More 
precisely, there exists $\Omega_\delta \in B(\Gamma)$ with 
$\Omega_\delta \ni id$ such that $\hat{d}(x,y) < \eta (\delta)$ as soon 
as $x_\gamma = y_\gamma$ for all $\gamma \in \Omega_\delta$. \\ 

Now fix $\Omega = \Omega_i$ and let $x, y \in M^\Gamma$ be such that 
$x_\gamma = y_\gamma$ for $\gamma \in \Omega \cdot \Omega_o$. Then
$$\bigl| f_o( \gamma^{-1} x ) - f_o( \gamma^{-1} y ) \bigr| < 
\delta$$ when $\gamma \in \Int_{D} (\Omega \cdot \Omega_o)$, where $D$ 
denotes the diameter of $\Omega_\delta$. Set $\hat{\Omega} = \Int_{D} 
(\Omega \cdot \Omega_o) \cap \Omega$ and decompose $f_\Omega$ into a convex linear 
combination as follows~: 

$$f_\Omega = \tfrac{|\hat{\Omega}|}{|\Omega|} f_{\hat{\Omega}} + 
\tfrac{|\Omega - \hat{\Omega}|}{|\Omega|} f_{\Omega - \hat{\Omega}}.$$
Then
$$\begin{array}{rclll}
\bigl| f_\Omega(x) - f_\Omega(y) \bigr| & \leq & \frac{ 
\hat{|\Omega|}}{|\Omega|} \bigl| f_{\hat{\Omega}} (x) - f_{\hat{\Omega}} (y) 
\bigr| &  + & \frac{|\Omega - \hat{\Omega}|}{|\Omega|} \bigl| f_{\Omega - 
\hat{\Omega}} (x) - f_{\Omega - \hat{\Omega}} (y) \bigr| \\
& \leq & \frac{ |\hat{\Omega}|}{|\Omega|} \;\; \delta & + &
\frac{|\Omega - \hat{\Omega}|}{|\Omega|} \\
& \leq & 2 \delta, 
\end{array}$$
provided the index $i$ is sufficiently large. Indeed, $\Omega - \hat{\Omega} \subset 
\partial_{D} \Omega = \partial_{D} \Omega_i$.  
\cqfd

\subsection{Repartition of classes according to degree and support}

Now we are ready to prove the nontriviality of $b$. It follows from \lref{PD} 
and \lref{delta} and does not further use the assumption that $X$ is a product.\\

\noindent
{\bf Proof of \pref{pos}} Let $A = A(\Omega_o)$ as before and let $(\Omega_i)$ be 
an amenable sequence in $\Gamma$. Lemma \ref{PD} implies that for each 
$i$ and each $a \in A_{\Omega_i}$, there exists a $c_a^{\Omega_i} \in 
[0,1]$ such that 
$$\varphi^{\star}_{A, \Omega_i, c_a^{\Omega_i}, \delta_{\Omega_i}} (a) 
\neq 0.$$ 
Thus any $a$ in $A_{\Omega_i}$ contributes to $b_{A, \Omega_i, c, 
\delta_{\Omega_i}}^{\ell, \nu}$ for some $c$ and $\ell$. Using the pigeon-hole 
principle, in the spirit of \pref{MI}, we will find some $\ell$ and some $c$ for which 
exponentially many classes of degree around $\ell|\Omega_i|$ are supported 
around $f_{\Omega_i}^{-1}(c)$. \\  

The degree of $a$ is an integer number between $0$ and $m | \Omega_o \cdot 
\Omega_i | \leq m \omega_o |\Omega_i|$, where $m = \Dim{M}$ and 
$\omega_o = |\Omega_o|$. Fix $r \in {\mathbb N}_o$. Then for each $i$, 
there exists an interval $J_i = [\frac{s-1}{r}, \frac{s}{r}] \subset 
[0, m \omega_o]$ for which the rank of the space $A_{\Omega_i}^{J_i}$ of 
classes in $A_{\Omega_i}$ whose degree belongs to the interval $J_i 
|\Omega_i| = \bigl[\frac{(s-1)}{r} |\Omega_i|, \frac{s}{r} |\Omega_i| \bigr]$ 
satisfies $$\rk A_{\Omega_i}^{J_i} \geq \frac{1}{r m \omega_o} \rk 
A_{\Omega_i}.$$     

\begin{lem}\label{slicing} Fix $k \geq 1$ and  $i \geq 1$. Then there exist 
an interval $I_i = [\frac{j-1}{k}, \frac{j}{k}] \subset [0,1]$ and a 
subspace $\ol{A}_i \subset A^{J_i}_{\Omega_i}$ such that  
\begin{enumerate}
\item[-] $a \in \ol{A}_i \Rightarrow c_a^{\Omega_i} \in I_i$,
\item[-] $\rk \ol{A}_i \geq \dfrac{1}{k} \rk A^{J_i}_{\Omega_i}$.
\end{enumerate}
\end{lem}

\Pf If $I \subset [0,1]$ and if $a$ is a cohomology class in $H^\star(M^\Gamma; F)$, 
denote by $r_I(a)$ the restriction of $a$ to the open set 
$f_{\Omega_i}^{-1}(I)$. Let $[0,1] = \cup_{j = 1}^{k} I_j$, with $I_j = 
[\frac{j-1}{k}, \frac{j}{k}]$. We decompose the space
$A^{J_i}_{\Omega_i}$ into a direct sum $$A^{J_i}_{\Omega_i} = A_1 \oplus 
\ldots \oplus A_k$$ in such a way as to satisfy the following properties~:  

$$A_k \oplus \Bigl( \Ker r_{I_k} \cap A^{J_i}_{\Omega_i} \Bigr) = 
A^{J_i}_{\Omega_i},$$ 
and for $j = 1, \ldots , k-1$,
$$\Biggl\{ 
\begin{array}{l} 
A_j \subset \Bigl( \Ker r_{I_{j+1} \cup \ldots \cup I_k} \cap 
A^{J_i}_{\Omega_i} \Bigl) \\
A_j \oplus \Bigl( \Ker r_{I_{j} \cup \ldots \cup I_k} \cap 
A^{J_i}_{\Omega_i} \Bigr) = \Bigl( \Ker r_{I_{j+1} \cup \ldots 
\cup I_k} \cap A^{J_i}_{\Omega_i} \Bigr). 
\end{array}$$ 
Thus if $a \in A_j$ then $c_a^{\Omega_i} \in I_j$.  Now there exists a 
$j = j(i)$ such that $\rk A_{j(i)} \geq \frac{1}{k} \rk A^{J_i}_{\Omega_i}$. 
Let $\ol{A}_i = A_{j(i)}$.

\cqfd

\noindent
The collection of intervals $J_i$ and $I_i$ being finite, 
there exist 
\begin{enumerate}
\item[-] a subsequence of $(\Omega_i)$, denoted $(\Omega_i)$ 
as well,
\item[-] an interval $J = [\ell - \frac{1}{2r}, \ell + \frac{1}{2r}]$, 
with $\ell = \ell(r)$,
\item[-] an interval $I = [c - \frac{1}{2k}, c + \frac{1}{2k}]$, with $c = c(k)$, 
\end{enumerate}
such that $\delta_{\Omega_i} \leq \frac{1}{4k}$, $J_i = J$ and $I_i = I$ 
for all $i$. Moreover, for each $i$ and each $a \in \ol{A}_{i}$, 
$$\supp a \subset f_{\Omega_i}^{-1} (- \infty, c_a^{\Omega_i} + 
\delta_{\Omega_i}) \subset f_{\Omega_i}^{-1} (- \infty, c + \tfrac{1}{k}).$$ 
Furthermore, \lref{PD} provides a class $b \in A_{\Omega_i}$ such that $$\supp 
b \subset f_{\Omega_i}^{-1} (c - \tfrac{1}{k}, + \infty) \;\; \mbox{ and } 
\;\; a \wedge b \neq 0.$$ Thus the map $\varphi^{\ell, \frac{1}{2r}}_{A, 
\Omega_i, c, \frac{1}{k}}$ is injective on $\ol{A}_i$ for all $i$. Therefore,

\begin{align*}
b^{\ell, \frac{1}{2r}}_{A} (c, \tfrac{1}{k}) & \geq \lim_{i \to \infty} 
\frac{1}{|\Omega_i|} \ln \bigl(\rk \ol{A}_i \bigr) \\
& \geq \lim_{i \to \infty} \frac{1}{|\Omega_i|} \ln \Bigl( 
\frac{1}{k} \frac{1}{r m \omega_o} \rk A_{\Omega_i} \Bigr) \\
& = \lim_{i \to \infty} \frac{1}{|\Omega_i|} \ln \Bigl( \frac{1}{k} 
\frac{1}{r m \omega_o} \bigl(\rk H^\star(M;F) \bigr)^{\Omega_i \cdot \Omega_o} \Bigr) \\
& = \lim_{i \to \infty} \frac{| \Omega_i \cdot \Omega_o |}{|\Omega_i|} 
\ln \bigl( \rk H^\star(M;F) \bigr) \\
& \geq \ln \bigl( \rk H^\star(M;F) \bigr).
\end{align*}   

To conclude, observe that $\ell = \ell(r)$ and $c = c(k)$. Let 
$\nu, \delta > 0$. There exists a subsequence $\ell(r_s)$ of 
$\ell(r)$ (\rp $c(k_l)$ of $c(k)$) such that $\ell(r_s)$ (\rp 
$c(k_l)$) converges to some $\ell \in \R^+$ (\rp $c \in [0,1]$).
Since $b^{\ell, \nu}_{A, \Omega} (c, \delta)$ increases with the 
size of $(\ell - \nu, \ell + \nu)$ and that of $(c - \delta, c + \delta)$, 

$$b^{\ell, \nu}_{A} (c, \delta) \geq b^{\ell(r_s), \frac{1}{2r_s}}_{A} (c(k_l), 
\tfrac{1}{k_l}) \geq \ln \bigl( \rk H^\star(M;F) \bigr),$$ 
for $l$ and $s$ sufficiently large. Thus 
$$b^\ell (c) \geq \ln \bigl( \rk H^\star(M;F) \bigr)  > 0.$$

\cqfd

\begin{rmk}\label{alg}{\rm
At least for products, $b_{A}^{\ell}(c) = b^\ell(c)$ provided $A$ 
contains $H^\star(M;F)$. Indeed, let $A^o = p_{\gamma}^* H^\star(M;F)$, 
some $\gamma$ in $\Gamma$, and let $A$ be 
another finite-dimensional subalgebra of $\a$ containing $A^o$, 
necessarily contained in some subalgebra $A^1 = p_{\Omega_1}^* 
H^\star(M^{\Omega_1};F)$ with $\Omega_1 \ni \gamma$. Thus $A^o_\Omega \subset 
A_\Omega \subset A^1_\Omega$ and 
$$b_{A^o, \Omega}^{\ell, \nu} (c, \delta) \leq b_{A, \Omega}^{\ell, 
\nu}(c, \delta) \leq b_{A^o, \Omega}^{\ell, \nu} (c, \delta) + \bigl(\rk 
H^\star(M;F)\bigr)^{|\Omega \cdot \Omega_1 - \Omega|}.$$    
Thus 
\begin{multline*}
b_{A^o}^{\ell, \nu}(c, \delta) \leq b_A^{\ell, \nu}(c,\delta) \leq \max \Bigl\{
b_{A_o}^{\ell, \nu}(c, \delta)\; , \\
\displaystyle{\lim_{i\to \infty} \tfrac{|\Omega_i \cdot \Omega_1 - 
\Omega_i|}{|\Omega_i|} \ln \bigl(\rk H^\star(M;F) \bigr) \Bigr\}  
= b_{A^o}^{\ell, \nu}(c, \delta).}
\end{multline*}
}\end{rmk}

\begin{ex}{\rm Let $F_o:S^1 \to [0,1]$ be a Morse function with two 
non-degenerate critical points, $x_-$, the minimum, and $x_+$, the maximum. 
Let $f_o= F_o \circ p_{id}$. The entropy of $f_o$ can be computed 
explicitely (cf.~\cite{L}, \S A4). Indeed, let $\mu$ denote the fundamental 
class of $M = S^1$. Then any class in $H^\star(M^\Omega;F)$ is of the type 

$$\displaystyle{\sum_{\Omega' \subset \Omega} n_{\Omega'} 
\mu^{\Omega'}},$$
where $n_{\Omega'} \in F$ and where $\displaystyle{\mu^{\Omega'} = \wedge_{\gamma 
\in \Omega'} (p_{\gamma}^* \mu)}$. Moreover, for a monomial $a = \mu^{\Omega'}$,

$$\left\{\begin{array}{l}
c_a^\Omega = \dfrac{|\Omega - \Omega'|}{|\Omega|}, \\
\deg a = |\Omega'|.
\end{array}\right.$$
}\end{ex}
The first equality is a consequence of the fact that the fundamental 
class $\mu$ can be supported in an arbitrarily small neighborood of any 
point (e.g.~$x_-$), implying that the class $\mu^{\Omega'}$ can be supported 
in any neighborhood of $M^{\Gamma - \Omega'} \times \{x_-\}^{\Omega'}$. 
It is now easy to convince oneself that if $A = p_{id}^* H^\star(M;F)$, then

$$\bi{\ell, \nu} = \# \bigl\{\Omega' \subset \Omega ; 
\left\{ \begin{array}{l} 
\dfrac{|\Omega - \Omega'|}{|\Omega|} \in (c - \delta, c + \delta), \\
\dfrac{|\Omega'|}{|\Omega|} \in (\ell - \nu, \ell + \nu). \\
\end{array}\right. \bigr\}.$$
Let $n = |\Omega|$. Then

$$\begin{array}{rcl}
\bi{\ell, \nu} & = & \displaystyle{\sum_{\substack{(c-\delta)n < j < 
(c+\delta)n \\ (1 - \ell - \nu)n < j < (1 - \ell + 
\nu)n}} \binom{n}{j}}.
\end{array}$$
Therefore

$$\begin{array}{rcl}
\displaystyle{\sup_{\substack{(c-\delta)n < j < (c+\delta)n \\ 
(1 - \ell - \nu)n < j < (1 - \ell + \nu)n}} 
\binom{n}{j} \quad \leq \quad \bi{\ell, \nu} \quad \leq \quad 2 \delta n
\sup_{\substack{(c-\delta)n < j < (c+\delta)n \\ 
(1 - \ell - \nu)n < j < (1 - \ell + \nu)n}} \binom{n}{j}}.   
\end{array}$$
Besides, by Stirling's formula,

$$\begin{array}{rcl}
\displaystyle{\frac{1}{n} \ln \binom{n}{j} } & \sim & 
\displaystyle{\frac{1}{n} \ln \Bigl(\frac{n^{n + \frac{1}{2}}}{j^{j + 
\frac{1}{2}} \; (n-j)^{n-j+\frac{1}{2}}} \Bigr) } \\ 
& \sim & \displaystyle{ - \bigl(\frac{j}{n}\bigr) \ln 
\bigl(\frac{j}{n}\bigr) -  \bigl(1 - \frac{j}{n}\bigr) \ln \bigl(1 - 
\frac{j}{n}\bigr) }.  
\end{array}$$
where we have removed terms that would produce a nul contribution in the 
limit. Now 

$$b_A^{\ell, \nu}(c, \delta) =   
\sup_{\substack{c-\delta < x < c+\delta \\ 1-\ell-\nu < x <
1-\ell+\nu}} \Bigl( -x \ln x - (1-x) \ln(1-x) \Bigr).$$
And thus  (using~\rref{alg})

$$b^\ell(c) = b_A^\ell(c) = \left\{ \begin{array}{lcl}
- \infty & \mbox{if} & c \neq 1 - \ell \\
- c \ln c - (1-c) \ln (1-c) & \mbox{if} & c = 1 - \ell.
\end{array}\right.$$
So $b$ is concentrated along the diagonal $c = 1-\ell$ and vanishes 
at the corners $(0,1)$ and $(1,0)$. The sum of the Betti number entropy 
is therefore given by 
$$b(c) = - c \ln c - (1-c) \ln (1-c).$$

\begin{rmk}{\rm
As suggested by the previous example, it is always true in the product 
case that, provided the functions $f_\Omega$ have constant range, the 
function $b(c)$ is nonnegative (i.e.~does not achieve the value $-\infty$). 
This is a consequence of the presence of the fundamental class whose 
support can be concentrated around any given point in $M$. 
}\end{rmk}

\section{Generalized Poincar\'e duality}

It has been observed in the product case that a class $a$ in $A_\Omega$, 
whose restriction to an open subset $O$ does not vanish, admits a nontrivial 
pairing with a class $b$ in $A_\Omega$ provided "$O$ is not too small", 
meaning is of type $p_\Omega^{-1} (O_o)$ for some $O_o$ in $M^\Omega$ 
(in the case $A = p_{id}^*H^\star(M;F)$). This suggests that  a condition 
generalizing Poincar\'e duality (more precisely its \lref{PD} version) should 
involve a filtration $(\T_{\Omega})_{\Omega \in B(F)}$ of the topology $\T$ of $X$ 
such that Poincar\'e duality holds in $A_\Omega$ for open subsets of 
$\T_\Omega$ (a precise statement follows). In the product case the topology 
$\T_\Omega$ is the one generated by the supports of the classes belonging to 
$A_\Omega$. It seems necessary to assume in addition that these topologies 
are induced by a family of (continuous) maps $(X, \T_\Omega) \to (X, \T)$ converging 
uniformely to the identity map. In the product case, for $A = p_{id}^*H^\star(M;F)$, 
these maps are the compositions $s_\Omega^o \circ p_\Omega$, where 
$s_\Omega^o$ is a section $M^\Omega \to M^\Gamma$ associated to a point 
$o \in M^\Gamma$ as follows~: 
$$\bigl( s_\Omega^o (x) \bigl)_\gamma = \left\{ 
\begin{array}{ccc}
x_\gamma & \mbox{if} & \gamma \in \Omega\\
o_\gamma & \mbox{if} & \gamma \notin \Omega.
\end{array}\right.$$
A last detail~: \lref{PD} holds when $A$ is the full cohomology algebra of a finite 
product. If $A$ is not of this type (e.g.~$M = {\mathbb T}^2$ and $A$ is generated 
by one cohomology class in $H^1({\mathbb T}^2)$), then the class $b$ does 
belong to $B_\Omega$ instead of $A_\Omega$, where $B = H^\star(M^{\tilde{\Omega}};F)$ 
and $\tilde{\Omega}$ is the smallest set for which  $H^\star(M^{\tilde{\Omega}};F) \supset A$. \\

Before stating the condition, we introduce the convention that whenever a sequence 
$something_\Omega$ (thus indexed by the set $B(\Gamma)$) is said to converge 
to $something$, it means that $something_{\Omega_i}$ converges to $something$ 
whenever $(\Omega_i)$ is an amenable sequence in $\Gamma$.
\begin{cond}\label{gPDc}
For any finite-dimensional subalgebra $A \subset \a$, there exist another finite-dimensional 
subalgebra $B$ with $A \subset B \subset \a$ and  a family of continuous maps 
$(r_\Omega : X \to X)_{\Omega \in B(\Gamma)}$ such that if $\T_\Omega$ denote the 
topology obtained by pulling back that of $X$ via $r_\Omega$, then

\begin{enumerate}
\item[-] $\gamma \circ r_\Omega = r_{\gamma \cdot \Omega} \circ \gamma$,
\item[-] $r_\Omega$ converges uniformly and monotonously to the identity map,
\item[-] for any $a \in A_\Omega$ and $O \in \T_{\Omega}$ such that $a|_{O} \neq 0$,
there exists a class $b \in B_\Omega$ with 
$$\left\{ 
\begin{array}{l}
\supp b \subset O \\ 
a \wedge b \neq 0.
\end{array}\right.$$
\end{enumerate}
\end{cond}
When, in addition to \aref{ass}, the previous condition is fulfilled, one may 
carry through the proofs of \lref{PD} and of \lref{delta}, and hence that of 
\pref{pos}. \\
 
Let $A$ be a finite-dimensional subalgebra of $\a$. Let $\delta>0$. Define 
$g_\Omega = f_\Omega \circ r_\Omega : (X, \T_\Omega) \to [0,1]$ and let 
$$\delta_\Omega = 4 \sup \{ | f_\Omega(x) - g_\Omega(x) | ; x \in X\}.$$

\begin{lem}\label{gPD}(Poincar\'e duality under \cdref{gPDc}) If $a$ belongs to 
$A_\Omega$ for some $\Omega \in B(\Gamma)$, then there exist a level 
$c_a^{\Omega}$ and a class $b$ in $B_\Omega$ such that  

\begin{enumerate}
\item[-] $\supp a \subset f_\Omega^{-1}(- \infty , c_a^{\Omega} + \delta_\Omega)$, 
\item[-] $\supp b \subset f_\Omega^{-1}(c_a^{\Omega} - \delta_\Omega, + \infty)$,   
\item[-] $a \wedge b \neq 0$.
\end{enumerate}
\end{lem}

\Pf As in the proof of \lref{PD} let 
$$c_a^{\Omega} = \inf \{ c \in [0,1] ; \supp a \subset f_\Omega^{-1}(-\infty, 
c) \}.$$
Now define $O_1 = f_\Omega^{-1} (c_a^\Omega - \delta_\Omega, + \infty)$ 
and $O_2 = f_\Omega^{-1} (c_a^\Omega -  \frac{1}{2}\delta_\Omega, + \infty)$. 
By construction,
$$\sup_{x \in X} \Bigl| f_\Omega (x) - g_\Omega (x) \Bigr| = \frac{\delta_\Omega}{4}.$$
Hence

$$O_2 \subset \bigl(g_\Omega \bigr)^{-1} (c - \tfrac{3}{4} \delta_\Omega, + \infty) 
\subset O_1.$$
Let $O$ denote $(g_\Omega)^{-1} (c - \frac{3}{4} \delta_\Omega, + \infty)$. 
Then $a|_{O} \neq 0$. Since $O \in \T_\Omega$, there exists a class $b \in 
B_\Omega$ with $\supp b \subset O$ and  $a \wedge b \neq 0$.

\cqfd

\begin{lem}\label{gdelta} The sequence $\delta_\Omega$ converges to $0$. 
In other words,  the sequence of functions $g_{\Omega} - f_{\Omega}$ converges 
uniformly to $0$. 
\end{lem}

\Pf 
Let $\delta > 0$. Since $r_\Omega$ converges uniformly and monotonously  
to the identity map and since $X$ is compact, there exists a finite set 
$\Omega_o \subset \Gamma$ containing $id$ such that for any $\Omega \supset 
\Omega_o$ for which the amenability ratio $\alpha(\Omega, \Omega_o)$ is 
sufficiently small,

$$\Bigl| f(r_{\Omega}(x)) - f(x) \Bigr| < \delta \qquad \forall x \in X.$$
Let $D_o = \diam \Omega_o$. If $\gamma \in \Int_{D_o} \Omega$, then

$$\Bigl| f(\gamma^{-1} r_\Omega(x)) - f(\gamma^{-1} x)\Bigr| =
\Bigl| f(r_{\gamma^{-1} \Omega}(\gamma^{-1} x)) - f(\gamma^{-1}  
x)\Bigr|  < \delta.$$
Let $\widetilde{\Omega} = \Int_{D_o} \Omega$. Then

$$\Bigl| f_{\widetilde{\Omega}}(r_\Omega(x)) - f_{\widetilde{\Omega}}  (x) \Bigr| 
\leq \frac{1}{|\widetilde{\Omega}|} \displaystyle{\sum_{\gamma \in \widetilde{\Omega}}} 
\Bigl|  f(\gamma^{-1} r_\Omega(x)) - f(\gamma^{-1} x) \Bigr| < \delta.$$
Hence 

$$\begin{array}{rl}
\Bigl| g_\Omega (x) - f_\Omega (x) \Bigr|  & \leq 
\displaystyle{\frac{|\widetilde{\Omega}|}{|\Omega|} \Bigl| f_{\widetilde{\Omega}} 
(r_\Omega (x)) - f_{\widetilde{\Omega}} (x) \Bigr|} \\
& \qquad \qquad + \displaystyle{\frac{|\Omega - \widetilde{\Omega}|}{|\Omega|} 
\Bigl| f_{\Omega  - \widetilde{\Omega}} (r_{\Omega} (x)) - f_{\Omega  - 
\widetilde{\Omega}}(x) \Bigr|} \\
& \leq \displaystyle{\frac{|\widetilde{\Omega}|}{|\Omega|} \delta + \frac{|\Omega - 
\widetilde{\Omega}|}{|\Omega|}} \\
& \leq 2 \delta,
\end{array}$$
Where the very last equality holds when, once more, $\Omega$ and $\Omega_o$ 
have a sufficiently small amenability ratio. 

\cqfd

Combining \lref{gPD} and \lref{gdelta} we obtain the following result, whose proof 
is essentially the same as that of \pref{pos}.

\begin{prop} The function $b$ achieves a strictly positive value.
\end{prop}

\section{Non-product example}\label{NPE}

Let $M$ be a projective algebraic variety, e.g.~the projective space ${\mathbb C}P^n$, 
and let $Y$ be a symbolic algebraic subvariety in $X = M^\Gamma$ (in the sense of 
\cite{G-TIDS} and \cite{G-ESAV}), that is, a compact subset $Y$ such that $Y_\Omega$, 
the "restriction" of $Y$ to each $\Omega$ in $B(\Gamma)$, defined as the image of the 
natural projection $M^\Gamma \to M^\Omega$, is an algebraic subvariety in $M^\Omega$. 
So $Y$ comes as the projective limit of the $Y_\Omega$'s, where one may (or may not) 
assume that $Y$ is $\Gamma$-invariant. We assume that for large enough 
$d(\Omega, \Omega')$ (depending on $Y$), the projection $Y_\Omega \cup Y_{\Omega'} 
\to Y_{\Omega \cup \Omega'}$ is onto. Observe that surjective maps between projective 
(in general K\"ahler) varieties are injective on the top-dimensional cohomology with complex 
coefficients due to existence of multivalued algebraic sections. (In fact if the fibers of such a 
map have, in a suitable sense, degree d, the same injectivity holds for $F_p$-coefficients, 
provided $p$ does not divide $d$). Therefore, if the target variety is non-singular, then the 
map is injective on all cohomology by Poincar\'e duality.
 
\begin{subex}{\rm
Let $M = {\C}P^n$ and consider a hypersurface in $M \times M$ represented by an equation 
$h(x,x') = 0$. Then the infinite chain of equations $h(x_i,x_{i+1}) = 0, i = ....,-1,0,1,...$  defines 
a subvariety $Y$ in $X = M^{\Z}$ invariant under the $\Z$-action.  
}\end{subex}
 
\begin{rmk}{\rm Unfortunately, even for generic $h$, it is unclear whether this $Y$ is non-singular 
in the sense that the restrictions of $Y$ to the intervals $[i,i+1,...,i+k]$, denoted $Y_{[i,i+1,...,i+k]}$, 
are non-singular.  However, a small (but non-$\Z$-invariant) perturbation $Y'$ of such a $Y$, 
allowing different $h$'s, i.e.~equations $h_i(x_i, x_{i+1}) = 0$, is non-singular by a simple argument 
(see \cite{G-TIDS} and \cite{G-ESAV}). Furthermore, the non-singular pertubations of $Y$ are all 
canonically homeomorphic and thus their cohomology can be attributed to $Y$ (alternatively, one 
may speak of a random $Y$ in $X$ with a suitable $\Z$-invariant probability measure on the space 
of strings $\{h_i\}$ and similarly introduce random potentials on $Y$ (and/or on $X$ itself). This 
significantly adds to possible examples and needs only a minor modification of our setting (with a 
reference to the sub-additive ergodic theorem).
}\end{rmk} 
 
\noindent
{\bf Continuation of the example.} The cohomology of our (desingularized) $Y$ enjoys the above 
product-like action on cohomology. In particular, for $\Gamma = \Z$, the homological entropy of 
a function exists. 

\begin{rmk}{\rm It seems hard to compute the (co)homologies of the above $Y_{[i,i+1,...,i+k]}$  or 
even to elucidate the properties of (the analytic continuation of) their entropic limit. However, it is 
easy to calculate the Chern numbers and thus the Euler characterictics of the Dolbeaut (and thus 
the ordinary) cohomology of all (desingularized!) $Y_{[i,i+1,...,i+k]}$.
}\end{rmk}

\section{Poincar\'e polynomial}

This section consists of defining the entropic Poincar\'e polynomial of a $\Gamma$-space.
It does not require the action to be product-like nor the group to be tileable. Amenability is the
only condition required here.\\

\noindent 
Consider the $A_\Omega$-Poincar\'e polynomial of $X$~: 

$$p_{A, \Omega}(t) = \sum_{d=1}^{\infty} t^d \rk A_\Omega^d,$$
where $A^d_\Omega$ denote the set of classes in $A_\Omega$ of (exact) 
degree $d$.

\begin{lem} The Poincar\'e polynomial of $X$ is $\Gamma$-invariant and subadditive, 
that~is
\begin{equation}\label{PP}
p_{A, \Omega \cup \Omega'} (t) \; \leq \; p_{A, \Omega}(t) \; p_{A, \Omega'}(t) \qquad 
\mbox{for} \qquad t \geq 0.
\end{equation}
\end{lem}

\Pf First observe that for any $\Omega_1, \Omega_2 \in B(\Gamma)$, the 
map

$$\displaystyle{\bigoplus_{d_1 + d_2 = d} A_{\Omega_1}^{d_1} \otimes 
A_{\Omega_2}^{d_2} \to A_{\Omega_1 \cup \Omega_2}^{d}}$$
is surjective. Thus 

$$\rk A_{\Omega_1 \cup \Omega_2}^d \leq \displaystyle{\sum_{d_1 + d_2 = d} 
\rk A_{\Omega_1}^{d_1} \rk A_{\Omega_2}^{d_2}},$$
which immediately implies the relation (\ref{PP}).

\cqfd

\noindent
Thus the limit 
$$\lim_{i \to \infty} \frac{1}{|\Omega_i|} \ln \bigl(p_{A, \Omega_i}(t)\bigr)$$
exists whenever $(\Omega_i)$ is an amenable sequence in $\Gamma$ (cf.~\cite{O-W}). We 
define the Poincar\'e polynomial of the group action $\rho : \Gamma \times M 
\to M$ to be 

$$p(t) = \sup_A \lim_{i \to \infty} \frac{1}{|\Omega_i|} \ln \bigl(p_{A, 
\Omega_i}(t)\bigr).$$  

\begin{rmk}{\rm
This definition is analoguous to that in \cite{G-TIDS} \S 1.14. 
Indeed, the process of factoring away $\eps$-fillable classes 
corresponds roughly to restricting to classes in $A_\Omega$.
}\end{rmk}


\begin{thebibliography}{xx}

\bibitem{G-TIDS} Misha Gromov, Topological invariants of dynamical systems 
and spaces of holomorphic maps. {\sl I. Math. Phys. Anal. Geom.} {\bf 
2} (1999), no.4, 323--415. 

\bibitem{G-ESAV} Misha Gromov, Endomorphisms of symbolic algebraic varieties. 
{\sl J. Eur. Math. Soc.} {\bf 1} (1999), no. 2, 109--197.

\bibitem{O-W} Donald S. Ornstein, Benjamin Weiss, Entropy and isomorphism 
theorems for actions of amenable groups. {\sl J. Analyse Math.} {\bf 48} (1987),
1--141.  

\bibitem{L} Oscar Lanford, Entropy and Equilibrium States in Classical 
Statistical Mechanics. {\sl Statistical mechanics and mathematical 
problems. Battelle Rencontres, Seattle, Wash.,} 1971. {\sl Lecture Notes 
in Physics,} {\bf 20}. Springer-Verlag, Berlin-New York, 1973. 

\bibitem{W} Weiss, Benjamin, Monotileable amenable groups. 
Topology, ergodic theory, real algebraic geometry, 257--262,
{\sl Amer. Math. Soc. Transl. Ser. 2}, {\bf 202}, AMS, Providence, RI, 2001. 

\end{thebibliography}
\end{document}